\documentclass{article}

\usepackage{amsmath, amsfonts,amssymb, amsbsy}
\newtheorem{theo}{\bf Theorem}[section]
\newtheorem{cor}{\bf Corollary}[section]
\newtheorem{defi}{\bf Definition}[section]
\newtheorem{nota}{\bf Notation}[section]
\newtheorem{lem}{\bf Lemma}[section]

\begin{document}
\begin{center}
{\bf {\Large

Nearly resolution V plans on blocks of small size.}}

\vskip10pt
 {\bf {\large Sunanda Bagchi\\
   Theoretical Statistics and Mathematics Unit\\
Indian Statistical Institute\\ Bangalore 560059, India.}}
\end{center}
\vskip10pt

\begin{abstract}
In  Bagchi (2010) main effect plans ``orthogonal through the block
factor" (POTB) have been constructed. The main advantages of a POTB
are that (a) it may exist in a set up where an ``usual" orthogonal
main effect plan (OMEP) cannot exist and (b) the data analysis is
nearly as simple as an OMEP.

   In the present paper we extend this idea and define the concept of
   orthogonality between a pair of factorial effects ( main effects
   or interactions)``through the block factor" in the context of a
   symmetrical experiment. We consider plans generated from an
   initial plan by adding runs. For such a plan
    we have derived necessary and sufficient conditions for a pair of
    effects to be orthogonal through the block factor in terms of
    the generators. We have also derived a sufficient condition on the
    generators so as to turn a pair of  effects aliased in the initial plan
    separated in the final plan.

       The theory developed is illustrated with plans for  experiments
       with three-level factors
    in situations where interactions between three or more factors are absent.
    We have constructed plans with blocks of size four and fewer runs
    than a resolution $V$ plan estimating all main effects and all but
    at most one two-factor interactions.

\end{abstract}


\section{Introduction}

A situation in which a treatment factor is neither orthogonal nor
confounded to a nuisance factor in a  main effect plan (MEP) was
first explored in Morgan and Uddin (1996). They considered a nested
row-column set up and derived a sufficient condition for a treatment
factor, possibly non-orthogonal to the nuisance factors, to be
orthogonal to every other treatment factor. They also derived
sufficient condition for optimality and constructed several series
of orthogonal and optimal MEPs.

     Mukherjee, Dey and Chatterjee (2001) Constructed main effect plans
  on blocks of a small size. Their plans also satisfy the
  property that  treatment factors are possibly non-orthogonal to
  the block factor and possess optimality property. However, their
  method relies on an existing orthogonal main effect plan (OMEP)
  as a starting point. Bose and  Bagchi (2007) provided examples
  of OMEPs in blocks of size 2 satisfying the same property as the
  plans of  Mukherjee, Dey and Chatterjee (2001) but
requiring fewer blocks. Extending this idea to an arbitrary block
size, Bagchi (2010) defined orthogonality between  a pair of
  treatment factors``through the block factor" and derived a
sufficient condition for that. This condition is weaker (thus
satisfied by more plans) than that of Mukherjee, Dey and Chatterjee
(2001). Method  of construction of plans orthogonal through the
block factor or ``POTB" were also presented, where it was seen that
a POTB may exist in a set up where an ``usual OMEP" cannot exist.

       In this paper we extend the concept of orthogonality through
the block factor to any pair of factorial effects (main effect or
interaction) in a symmetric experiment. We compare this concept with
the usual orthogonality. We also obtain necessary and sufficient
condition on a plan which makes the inference on a given factorial
effect free from the involvement of all other factorial effects, in
the presence or absence of a blocking factor. [See Theorem \ref
{cond-orth} and Remark 3.2]

Next we concentrate on construction. We consider  plans obtained by
``expanding an initial plan along a subspace" [see Definition \ref
{gen.plan}]. We derive sufficient conditions on the subspace so that
the relation between the effects is improved upon in the final plan.
That is, a pair of effects aliased in the initial plan is not so in
the final plan. Similarly, a pair aliased or non- orthogonal in the
initial plan becomes orthogonal (through the block factor) in the
final plan.

We have  illustrated these ideas using plans for three-level factors
assuming all interactions involving three or more factors to be
absent. We have obtained plans on blocks of size four for $3^3$
experiment on 6 blocks, $3^4,3^5$  experiments, each on 18 blocks
and $3^6$ experiment on 24 blocks. The first two plans estimate all
effects assumed in the model, while the last two plans estimate all
but one two-factor interaction. [See Theorems \ref {3^3}, \ref
{3^4}, \ref {3^5} and \ref {3^6}]. Each plan requires fewer runs
than a resolution V plan for the corresponding experiment.

\section{Preliminaries}

Let $s$ be a prime power. Let us recall the terminology of the
$m$-dimensional Euclidean geometry and a few other terms and
notations required for a symmetric experiment. We shall mostly
follow the notations of Bose (1947).

\begin{nota} \label{EGns}
(a) $F$ will denote the Galois field of order $s$. $0$ will denote
the additive identity of $F$.


(b) $F^m$ will denote the vector space of dimension $m$ over $F$. We
shall think of the vectors in $F^m$ as column vectors.

 (c) The points of an $m$-dimensional
 Euclidean geometry ($EG(m,s)$) are the vectors of $F^m$. A point
 is denoted by $x = (x_1, \cdots x_m)', \: x_i \in F, \: i=1, \cdots
 m$.

 (b) A hyperplane is a coset of an $(m-1)$-dimensional subspace of
$F^m$. Specifically, for $a \in F^m,\:a \neq 0$ and $t \in F$, the
hyperplane $H_a(t)$ is the set of points of $EG(m,s)$ satisfying
$a'x = t$.

(c) A pencil is a set of $s$ parallel hyperplanes of $EG(m,s)$. The
pencil $P_a$ denotes the set of parallel hyperplanes $\{H_a(t), t
\in F\}$. Thus, the pencils are the sets $P_a = \{H_a(t), t \in F\},
\: a \in F^m \setminus\{0\}$.
\end{nota}

{\bf Remark 2.1 :}  The pencil  $P_a$ is the same as the pencil
$P_b$, if $b = pa, p \in F, p \neq 0$.

\begin{nota} \label{expt}
Consider an $s^m$ experiment on $n$ runs. Let us name the factors as
$A_i, i=1,2, \cdots m$.

(i)  A point $x =(x_1,x_2, \cdots x_m)$ in $EG(m,s)$
 represents a level combination (run),  in which $A_i$ is at level $x_i,
i=1,2, \cdots m$.

(ii) The combined effect of all the main effects and interactions
present on the level combination $x$ will be denoted by $\tau_x$.

(ii) By a (factorial) effect we mean a main effect or an interaction
and it will generally be denoted by $D,E$ etc.

(iii) For a vector $a =(a_1, \cdots a_m)' \in F^m$, $ E_a$ will
denote the factorial effect $A^{a_1}_1A^{a_2}_2 \cdots A^{a_m}_m$.
Thus, the main effect of $A_i$ will be denoted by $ E_a$, where $a$
has 1 in only the $i$th position and $0$ elsewhere.
 Similarly, for a $3^5$ experiment the interaction $A_2 A^2_4$ will
be denoted by $E_b$, where $b = (0,1,0,2,0)'$.


(iv) The pencil $P_a$ will represent the factorial effect $ E_a$.

(v) We shall say that a run $x$ is in  the $t$th level of the
factorial effect $ E_a$ if $x$ satisfies $a'x =t$.
\end{nota}


{\bf Remark 2.2:} Recall that a contrast belonging to the factorial
effect $E_a$ is of the form $ \sum\limits_{t \in F} l_t
\sum\limits_{x \in H_a(t)} \tau_x$, where  $l_t$ is a real number
for each $t \in F$, such that $\sum\limits_{t \in F} l_t =0$. Thus,
for $b = pa, p \in F, p \neq 0$, a contrast belonging to the effect
$E_a$ also belongs to $E_b$. This is consistent with the property of
pencils noted in Remark 2.1.

\begin{nota} \label{incidence} Consider an $s^m$ experiment on $b$
blocks of size $k$ each.

(i) ${\cal E}$ will denote the set of effects believed to be
present.

(ii) Consider an effect $E \in {\cal E}$. The replication number
$r^E_t$ of the level $t$ of $E$ is the number of runs in which $E$
is at level $t$. The replication vector $r^E$ is the $s \times 1$
vector with $r^E_t$ as the $t$th entry, $t \in F_s$.

(iii) For two effects $D$ and $E$, the $D$-versus $E$ incidence
matrix
 is a $s \times s$ matrix denoted by $N^{DE} = ((n^{DE}_{pq}))_{p,q \in F}$,
 where $n^{DE}_{pq}$ denotes the number of runs which are in the $p$th level of $D$ as
well as the $q$th level of $E$.

Note that when $D$ and $E$ are main effects, say of factors $A_i$
and $A_j$, $n^{DE}_{pq}$ is the number of runs in which the $i$th
entry is $p$ and the $j$th entry is $q$.


(iv) For an effect $E$, the $E$-versus block incidence matrix
 is a $ s \times b$ matrix denoted by $N^{E,bl} = ((n^{E,bl}_{pj}))_{p \in F, 1\leq j \leq b}$,
 where $n^{E,bl}_{pj}$ denotes the number of runs which are
 in the $p$-th level of $E$ and in the $j$th block.
 \end{nota}

 We now define the concept of {\bf  orthogonality through the block
factor}.

\begin{defi} \label{orthbl} Two effects $D$ and $E$ are said to be
 {\bf  orthogonal through the block factor (OTB)} if
\begin{equation}  \label{orthbl-inc}
k {\mathbf N}^{DE} = {\mathbf N}^{D,bl} ({\mathbf N}^{E,bl})' .
\end{equation}

We denote this by $D \bot E (bl)$.
\end{defi}

{\bf Remark 2.3:} Condition (\ref {orthbl-inc}) is equivalent to
equation (7) of Morgan and Uddin (1996) in the context of nested
row-column designs.

The next result follows immediately from the definition of OTB.

\begin{theo}\label{condOrth1} A sufficient condition for
 effects $D$ and $E$  to be `OTB' is that each of the three incidence
 matrices ${\mathbf N}^{D,bl}, {\mathbf N}^{E,bl}$ and  ${\mathbf
 N}^{DE}$ have all entries equal. Specifically, the condition is
 \begin{eqnarray*}{\mathbf N}^{D,bl} & = & (n/(bs) J,\\
{\mathbf N}^{E,bl} & = & (n/(bs) J\\
\mbox{ and }{\mathbf N}^{DE} & = & (n/s^2) J. \end{eqnarray*}

Here $J$ is the all-one matrix.
\end{theo}

Before closing this section we present a few well-known results of
block designs.

Consider an equireplicate block design $d$ with   $v$ treatments and
$b$ blocks of size $k$ each. Let the replication number be $r$.

\begin{theo}\label{N-repl} The incidence matrix $N$ of $d$
 satisfies $$N 1_b = r \mbox{ and } N'1_v = k 1_b.$$\end{theo}

 Consider a pair of  equireplicate block designs $d_1$ and $d_2$ on
 the same set of $v$ treatments on $b$ blocks of size $k$ each.
 Let their incidence matrices be denoted by $N_1$ and $N_2$
 respectively. Let the common replication number be $r$. Then,
 clearly we have :

 in view
 of Lemma \ref {N-repl} we can say the following.

\begin{theo}\label{pairInc} $$N_1 N'_2 1_v = r k 1_v.$$\end{theo}

\section{Orthogonality through block versus usual orthogonality}

We shall now compare OTB with the usual orthogonality. For the sake
of simplicity we assume a main effect plan, but the results are
valid for any factorial effect.

We first present a few notation and well-known results.

\begin{nota} \label{model} (a) Consider a blocked MEP for $m-1$ factors.
The block factor is named as $A_m$ and the general effect as
$A_{m+1}$. The model is expressed as
\begin{equation}\label{model} {\mathbf Y} =  {\mathbf X} \beta,
\mbox{ where } {\mathbf X}  = \left[ \begin{array}{ccc}  {\mathbf
X}_1 & \cdots &  {\mathbf X}_{m+1}\end{array} \right ]  \mbox{ and
}\beta =\left[ \begin{array}{c} \alpha^1\\ \vdots
\\\alpha^{m+1}\end{array} \right ].\end{equation}

Here, $X_i$, the {\bf design matrix} for $A_i$ is a $0-1$ matrix -
the $(u,t)$th entry of ${\mathbf X}_i $ is 1 if in the uth run the
factor $A_i$ is set at level t and 0 otherwise, $ i=1,2, \cdots m$
and $X_{m+1} = {\mathbf 1}_n$.

Throughout this section $r^i$  will denote the vector of
replications of $A_i$, $N_{ij}$ the $A_i$-versus $A_j$ incidence
matrix and $R_i$ will replace $N_{ii}$.

(b)  $ a_i \times 1$ vector $\alpha^i$ will denote the vector of
unknown effects of $A_i,\: 1\leq i \leq m+1$. Thus, $\alpha_{m}$ is
the vector of block effects and $\alpha_{m+1}$ is the general
effect.

(c)  For any $m \times n$ matrix $A$, ${\cal C} (A)$ will denote the
column space of $A$. Further,  $P_A$ will denote the projection
operator on the column space of $A$. In other words, $P_A = A
(A'A)^- A'$, where $B^-$ denote a g-inverse of $B$.

(d)  Let $I =  \{1,2, \cdots m +1\}$ and   $S = \{i,j, \cdots\}$ be
  a subset of $I$.  We shall use the following
 notation for the sake of compactness.

 (i) $ X_S$ will denote $ \left[ \begin{array}{ccc}{\mathbf X}_i &
 {\mathbf X}_j & \cdots \end{array} \right ]$.

(ii) $ \alpha^S$ will denote  $ \left[ \begin{array}{c} \alpha^i \\
\alpha^j \\ \vdots \end{array} \right ]$.

(iii)  $ P_i$ will denote the projection operator onto the column
space of $ X_i, i \in I$. Further, $ P_S$ will denote the projection
operator onto the column space of $ X_S$.

(e) Sum of squares : Fix a set of factors $T$ of $I$. For i not in
$T$, we define $SS_{i;T}$, the sum of squares for $F_i$, adjusted
for the factors $F_t, t \in T$.  More generally we define the
combined sum of squares for the set factors $\{F_i, i \in S\}$,
 $SS_{S;T}$,  adjusted for the factors $F_t, t
\in T$ (T disjoint from S) as follows.
\begin{eqnarray*} SS_{i;T} & = & Q'_{i;T} (C_{ii;T})^- Q_{i;T} \\
\mbox{ and }  SS_{S;T} & = &Q'_{S;T} (C_{S,S;T})^- Q_{S;T}
\end{eqnarray*}

\end{nota}
For ready reference, we present the following well-known results.

\begin{lem} \label{ssadj}
(a) For $ i \neq j,\: i,j=1, \cdots m+1$,
  $SS_{i;j}$ is the quadratic form $Y'P_U Y$, where $U
=(I-P_j)X_i$.

(b) More generally, for two disjoint subsets $S$ and $T$ of I, $
SS_{S;T}$ is the quadratic form $Y'P_V Y$, where $V =(I-P_T) X_S$.

(c) The so-called unadjusted sum of squares for $F_i$ is $SS_{i;m+1}
= T'_i (R_i)^{-1} T_i -G^2/n$. \end{lem}

\begin{lem} \label{PA-PD} Consider a matrix $W$ partitioned as
$\left[ \begin{array}{cc} U & V \end{array} \right ]$. Let $ Z = (I
- P_V)U$. Then, $P_W - P_V = P_Z$. \end{lem}

\begin{cor} \label{Ti} Let $ T \subset I$ and $i \in I \setminus T$.
 Let $D = (I - P_T)X_i$. Then,
$$P_D = P_{T*} - P_T, \mbox{ where } T* = T \cup \{i \}.$$
\end{cor}

\begin{nota} \label{sumsofsq}(a) The total sum of squares and the error sum of squares
will be denoted by $SS_{tot}$ and $SS_E$ respectively.

(b) Fix a factor $A_i, 1 \leq i \leq m$. Let $T = \{i+1, \cdots m+1
\}$ and $\bar{i} = I \setminus \{i \}$.

(i) Let $SS_{i;all>} = SS_{i;T}$ and

(ii) $SS_{i; all} = SS_{i;\bar{i} }$.

Thus, $SS_{i;all>}$ is the sum of squares for $A_i$, adjusted for
the factors $ A_{i+1}, \: \cdots A_{m+1}$, while $SS_{i; all}$
denotes the sum of squares for $A_i$, adjusted for all other
factors.
\end{nota}

Now we seek the answer to the following questions. {\bf Consider a
main effect plan for m factors ($m \geq 3$). Fix a factor, say
$A_i$. What conditions the design matrices must satisfy so that the
sum of squares for $A_i$ adjusted for all others is the same as

(a) the unadjusted sum of squares for $A_i$ ?

(b) the sum of squares for $A_i$ adjusted for only one factor, (say
$F_m$) ?}

[That is so far as $A_i$ is concerned, other factors are virtually
absent.]

\begin{theo} \label{cond-orth} Fix a factor, say $A_i$.

(a)A necessary and sufficient condition for $SS_{i;all} =
SS_{i;m+1}$ is that  the incidence matrix  $N_{ij}$ satisfies the
proportional frequency condition of Addelman (1962). That is $
N_{ij} = r_i r'_j/n $ , for every $j \neq i$ .

(b) A necessary and sufficient condition for $SS_{i;all} = SS_{i;m}$
is that
\begin{equation} \label{orth-inc} N_{ij} = N_{im} (R_m)^{-1} N'_{jm},
j \neq i,\; 1 \leq i,j \leq m-1. \end{equation}
\end{theo}

[Recall Notation \ref {model}]

The proof relies on a lemma that we present now.

\begin{lem} \label{P1A&P2A} Consider matrices $A (m \times n), B((m
\times p)$  such that
$$ {\cal C} (B) \subseteq {\cal C} (A).$$

Let $C((m \times q)$ be any matrix. Then the necessary and
sufficient condition that $ {\cal C} (P_B C) = {\cal C} (P_A C)$ is
that $ (P_A - P_B) C = 0$. \end{lem}

{\bf Proof of theorem \ref {cond-orth}:} Let $T = \{ 1,2, \cdots
i-1, i+1, \cdots m-1\}$,  $T^* =  T \cup \{m \}$ and $T^{**} = T^*
\cup \{m+1 \}$. From Lemma \ref {ssadj}
 and Notation \ref {sumsofsq}, we see that
$$ SS_{i;all} = Y'P_U Y,\; SS_{i;m+1} = Y'P_V Y,\; SS_{i;m} =Y'P_W Y, $$

 where $U = (I - P_{T^{**}}) X_i, \;  V = (I - P_{m+1})X_i$ and $W = (I
 - P_m)X_i$.

{\bf Proof of (a);} From the expressions above, a necessary and
sufficient condition for $SS_{i;all} = SS_{i;m+1}$ is that $ P_U =
P_V$.

Therefore, in view of  lemma \ref {P1A&P2A} , the necessary and
sufficient condition for $SS_{i;all} = SS_{i;m+1}$ is that
\begin{equation}\label{cond} (P_{T^*} - P_{m+1}) X_i = 0. \end{equation}

Now by Lemma \ref {PA-PD}, $P_{T^*} - P_{m+1} = P_Z$, where $Z = (I
- P_{m+1}) X_T$. Thus, (\ref {cond}) is $ \Leftrightarrow P_z X_i =
0\;  \Leftrightarrow X'_i(I- P_{m+1})X_T =0 \;\Leftrightarrow X'_i
(I - P_{m+1}) X_j = 0,  j \neq i,$ which is the same as the
proportional frequency condition.

{\bf Proof of (b) :}  Proceeding along the lines as in the proof of
(a), we find that the necessary and sufficient condition for
$SS_{i;all} = SS_{i;m}$ is that
\begin{equation} P_Z X_i = 0,\; \mbox{ where } Z = (I - P_m) X_T.
\end{equation}

But this condition $\Leftrightarrow X'_i (I - P_m) X_T =0
\Leftrightarrow X'_i (I - P_m) X_j = 0, j \neq i,\; 1 \leq i,j \leq
m-1$.  This condition can be expressed in the same form as in the
statement by using the relation between the design matrices and
 the incidence matrices. $\Box$

{\bf Remark 3.1 :} The sufficiency part of (a) of Theorem \ref
{cond-orth} is already known. We have now shown that the condition
is also necessary.

{\bf Remark 3.2 :} The results of Theorem \ref {cond-orth} can be
easily generalized to a situation where interactions are present.
Part (a) says that the sums of squares for testing the significance
of  effect $E$ (say $SS_{E;adj}$) is the same as the so-called
unadjusted sum of squares for $E$ if and only if $E$ satisfies PFC
with every other effect present. Similarly, part (b) says that
$SS_{E;adj}$ is the same as that sums of squares for $E$ adjusted
for the block effect if and only if $E \bot D (bl)$ for every other
effect $D$.

\section{Construction}
  In the rest of this paper by orthogonality we mean orthogonality
   through the block factor.

    Suppose our aim is to construct a plan for an $s^m$ experiment on blocks
of size $k$ each. We begin with an initial plan ${\cal P}$ and then
develop or expand it with the help of a vector subspace $V$ of $F^m$
to our final plan. In this section we find out conditions on $V$, so
that the final plan satisfies certain desirable properties.

\begin{nota} Consider a plan ${\cal P}$ for an $s^m$ experiment on $b$ blocks
of size $k$ each. Let $n = bk$. Thus, there is a total of $n$ runs
in the plan, which in general may not be distinct.

(a) Let ${\cal B}$ denote the set of all blocks in ${\cal P}$, say
${\cal B} = \{ B_j, \: j=1, \cdots b\}$.

 (b) Let ${\cal E} = \{E_a, a \in I\}$, be the set of effects believed to
be present. Thus, $I \subseteq F^m \setminus \{0\}$.

(c) For $ a \in I$, the replication number $r^a_t$ of the level $t$
of effect $E_a$ is the number of runs in ${\cal P}$ in which $E_a$
is at level $t$. Thus, $r^a_t =  \sum_{j=1}^b \sum \limits_{x \in
B_j : a'x =t} 1$.

The replication vector $r^a$ is the $s \times 1$ vector with $r^a_t$
as the $t$th entry, $t \in F_s$.

(d) For $a,b \in I$, the $E_a$ versus $E_b$ incidence matrix is the
$ s \times s$ matrix $N^{ab}$, where $N^{ab} = ((n^{ab}_{\alpha,
\beta}))_{\alpha, \beta \in F}$, and $n^{ab}_{\alpha, \beta}$ is the
number of runs in ${\cal P}$ in which $E_a$ is at level $\alpha$ and
$E_b$ is at level $\beta$. That is

$$n^{ab}_{\alpha, \beta} =  \sum_{j=1}^b \sum
\limits_{\underset{a'x = \alpha, b'x = \beta } { x \in B_j :}  }
1.$$

(e) For $a \in I$, the $E_a$-versus block incidence matrix $L^a$ is
the $s \times b$ matrix $L^a = ((l^a_{\alpha, j}))_{\alpha \in F, 1
\leq j \leq b}$, where $l^a_{\alpha, j}$ is the number of runs
 in the $j$th block of ${\cal P}$ in which $E_a$ is at level
$\alpha$. That is

$$l^a_{\alpha, j} =  \sum \limits_{x \in B_j :  a'x = \alpha} 1.$$
\end{nota}

\begin{defi} \label{gen.plan} Let $V$ be a vector subspace of $F^m$,
say of dimension $t \geq 1$. For any block $B \in {\cal B}$, and any
vector $v \in V$, $v + B $ will denote the set $\{v + x : x \in B
\}$ of runs. By the {\bf expansion} $V({\cal P})$ of ${\cal P}$
along $V$ we shall mean the plan (for the same experiment, with
$b.s^t$ blocks of size $k$ each) whose set of blocks is $\{ v + B :
v \in V, B \in {\cal B} \}$.
\end{defi}

\begin{nota} We shall use the following notation for the replication
vectors and incidence matrices for  an expansion $V({\cal P})$ of a
plan ${\cal P}$ along a vector subspace $V$.

(a) $\tilde{r}^a$ will denote the replication vector of  $E_a$ in
$V({\cal P})$.

(b) $\tilde{N}^{ab} = ((\tilde{n}^{ab}_{\alpha,\beta}))_{\alpha,
\beta \in F} $ will denote the $E_a$ versus $E_b$ incidence matrix
in $V({\cal P})$.

(c) Let $J = \{1, 2, \cdots b \}$. The $s \times b|V|$ matrix
$\tilde{L}^a$ will denote the $E_a$ versus block incidence matrix in
$V({\cal P})$. Its entries will be denoted by
$\tilde{l}^a_{\alpha,\tilde{j}},\: \alpha \in F, \: \tilde{j} \in J
\times V$.
\end{nota}

We now find the replication vectors  and incidence matrices of
effects in $V({\cal P})$ in terms of those in the original plan
${\cal P}$.

\begin{theo} \label{IncV(P)}Consider an expansion of a plan
${\cal P}$ along a vector subspace $V$ as in Definition \ref
{gen.plan}. Consider a pair of effects $E_a$ and $E_b$. Then, the
replication vectors of these effects and the incidence matrices
involving $E_a,E_b$ and block for  $V({\cal P})$ are given in terms
of the corresponding quantities for ${\cal P}$ as follows.

\begin{eqnarray} \label{gen.Plan}
\tilde{n}^{ab}_{\alpha,\beta}  & = &\sum \limits_{ v\in V}
n^{ab}_{\alpha - a'v, \beta - b'v},\; \alpha, \beta \in F. \\
\tilde{r}^a_{\alpha} &= &\sum \limits_{v\in V} r^a_{\alpha - a'v},\; \alpha \in F.\\
\tilde{l}^a_{\alpha,\tilde{j}} &=& l^a_{\alpha - a'v,j}, \; \alpha
\in F, \: \tilde{j} = (j,v), \: 1 \leq j \leq b, \: v \in V.
\end{eqnarray}

\end{theo}

\begin{cor}\label{incbl} For a plan $V({\cal P})$ which is an expansion
of ${\cal P}$ along $V$, the following hold. For $\alpha, \beta \in
F$, the $(\alpha,\beta)$th element of $\tilde{L}^a (\tilde{L}^b)'$
is given by

$$\sum \limits_{ v\in V} (L^a (L^b)')_{\alpha - a'v,\beta - b'v}.$$
\end{cor}

 We now find the condition for orthogonality between a pair of effects in the final
 plan.

\begin{nota} For $V$ as in Definition \ref {gen.plan} let $ V^{\bot}$
denote its orthocomplement. That is $ V^{\bot}$ is vector subspace
of $F^m$ of dimension $m-t$, given by $ V^{\bot} = \{ w \in F^m :
w'v = 0, \forall v \in V\}$.
\end{nota}

\begin{defi}\label{effectGroup} Given a vector subspace $V$ of
$F^m$, let $\sim_V$ denote the binary relation on the set of all
effects given by

$$E_a \sim_V E_b \mbox{  if  } <a> + V^{\bot} = <b> + V^{\bot} .$$

Clearly $\sim_V$ is an equivalence relation and so it partitions the
set of all effects into the corresponding equivalence classes. We
define ``the effects classes relative to $V$" to be the
$\sim_V$-equivalence classes. \end{defi}

We now present our main result.

\begin{theo}\label{condOrthg} Consider an expansion $V({\cal P})$ of
${\cal P}$ along $V$. Fix a pair of effects $E_a$ and $E_b$.

(a) If $E_a$ and $E_b$ are from different effects classes relative
to $V$, then $E_a \bot E_b (bl)$ in $V({\cal P})$.

(b) Suppose $E_a$ and $E_b$ are from the same effects classes
relative to $V$. Then, the following hold.

(i) If both $a$ and $b$ are in $ V^{\bot}$, then the relation
between $E_a$ and $E_b$ in $V({\cal P})$ is the same as that in
${\cal P}$. That is  $E_a$ is confounded, non-orthogonal or
orthogonal to $E_b$ in $V({\cal P})$ if it is so in ${\cal P}$.

(ii) If neither $a$ nor $b$ is in $ V^{\bot}$ and $E_a$ is aliased
with $E_b$ in ${\cal P}$, then $E_a$ and $E_b$ are no more aliased
in $V({\cal P})$, provided $V$ is non-trivial.
\end{theo}

The proof is based on two lemmas we present now.

\begin{lem}\label{UseHypPlane} Let $V$ be a  $t$-dimensional
vector subspace of $F^m,\: a,b \in F^m$ and $\alpha, \beta \in F$.
Consider the subset $S(a,b, \alpha, \beta)$ of $V$, namely $S = \{ v
\in V : a'v = \alpha, \: b'v = \beta\}$. Then, the following hold.

(a) $S$ is empty in the following cases. (i) $a \in  V^{\bot}, \:
\alpha \neq 0$, (ii) $b \in  V^{\bot}, \:  \beta \neq 0$ and (iii)
$a - cb \in  V^{\bot}$, for some  $c \in F, c \neq 0$, but $\alpha
\neq c\beta$.

(b) $|S| = s^t$ if $a,b \in  V^{\bot}, \:  \alpha = \beta = 0$.

(c) $|S| =  s^{t-1}$ if one of the following conditions is
satisfied. (i) $ a \in  V^{\bot}, \: \alpha = 0, b \not\in
V^{\bot}$, (ii) $b \in  V^{\bot}, \:  \beta = 0, \: a \not\in
V^{\bot}$, (iii) $a \not\in  V^{\bot},\: b \not\in V^{\bot}$ and
$\exists \: c \in F$  such that $ a - cb \in V^{\bot}, \: \alpha =
c\beta$.

(d)  $|S| =  s^{t-2}$ if $a  \not\in  V^{\bot},\: b  \not\in
V^{\bot},\: <a> + V^{\bot} \neq <b> + V^{\bot}$ .\end{lem}

{\bf Proof :} This is trivial if $ a \in  V^{\bot}$ or $b \in
V^{\bot}$. So, we assume $a  \not\in  V^{\bot},\: b  \not\in
V^{\bot}$.

{\bf Case 1:}  $<a> + V^{\bot} = <b> + V^{\bot}$. Then, $\exists c
\in F$ such that $a - cb \in V^{\bot}$. If $ \alpha \neq c\beta$,
then the set is empty, proving (a)(iii). If $\alpha = c\beta$, then
the set is nothing but $\{ v \in V,\:a'v = \alpha\}$, which has size
$s^{t-1}$. Thus, (c) (iii) is proved.

{\bf Case 2 :} $<a> + V^{\bot} \neq <b> + V^{\bot}$. Then, $v
\longmapsto a'v$ and $v \longmapsto b'v$ are non-zero linear
functionals on on $V$. Their kernels are  $(t-1)-$dimensional vector
subspaces of $V$, namely $a^{\bot} \cap V$ and $b^{\bot} \cap V$.
These subspaces are distinct, since $<a> + V^{\bot} \neq <b> +
V^{\bot}$. The sets $\{ v \in V,\:a'v = \alpha\}$ and $\{ v \in V,\:
\: b'v = \beta\}$ are cosets of these subspaces and therefore are
euclidean hyperplanes in the $t$-dimensional euclidean space over
$F$. Since the vector subspaces are distinct, these hyperplanes are
distinct and non-parallel. But any two non-parallel hyperplanes in
the $t$-dimensional euclidean space meet in an euclidean subspace of
dimension $t-2$, having $s^{t-2}$ common points. This proves (d) and
hence completes the proof.$\Box$.

\begin{lem}\label{PrpIncMat} Let $a,b$ and $V$ be as in Lemma \ref
{UseHypPlane}. For any $s \times s$ matrix $M$ with rows and columns
indexed by $F$, we define the $s \times s$ matrix $\tilde{M}$ as
follows.

$$\tilde{m}_{\alpha, \beta} = \sum \limits_{v \in V}
 m_{\alpha - a'v,\beta - b'v}, \; \alpha, \beta \in F.$$
Then, we have :-

(a) If $a, b \in  V^{\bot}$, then, $\tilde{M} = s^t M$.

(b) If $a \in V^{\bot},\: b \not\in  V^{\bot}$, then,
$\tilde{m}_{\alpha, \beta} = s^{t-1}\sum \limits_{\gamma \in F}
m_{\alpha, \gamma}, \; \alpha, \beta \in F.$

(c)  If $a \not\in V^{\bot},\: b \in  V^{\bot}$, then,
$\tilde{m}_{\alpha, \beta} = s^{t-1}\sum \limits_{\gamma \in F}
m_{\gamma,\beta}, \; \alpha, \beta \in F.$

(d) If  $a \not\in V^{\bot},\: b \not\in  V^{\bot}, \:  <a> +
V^{\bot} = <b> + V^{\bot}$, then,  $\tilde{m}_{\alpha, \beta} =
s^{t-1}\sum \limits_{u  \in F} m_{\alpha - u, \beta - cu}$, for some
$c  \neq 0 \in F$.

(e) If $a \not\in V^{\bot},\: b \not\in  V^{\bot}, \:  <a> +
V^{\bot} \neq <b> + V^{\bot}$, then,  $\tilde{M} = c s^{t-2}J$,
where $c$ is the sum of all entries of $M$ and $J$ is the all-one
matrix.
\end{lem}

{\bf Proof :} From the definition of  $\tilde{M}$, its
$(\alpha,\beta)$th entry is nothing but

$$\sum \limits_{z,w \in F} m_{z,w}. \|\{v \in V,\: a'v = \alpha - z,
b'v = \beta - w \} \|. $$

Therefore the result is immediate from Lemma \ref
{UseHypPlane}.$\Box$

{\bf Proof of Theorem \ref {condOrthg}:} Put $P = L_a L'_b, M =
N^{ab}$. Then, by Theorem \ref {IncV(P)}and Corollary \ref  {incbl},
we have (in the notation of Lemma \ref {PrpIncMat})
 $\tilde{L_a} \tilde{L'_b}  = \tilde{P}$ and  $\tilde{N^{ab}} =  \tilde{M}$.

Proof of (a): Consider two cases.

{\bf Case 1: }  Exactly one of $a$ and $b$ is in $ V^{\bot}$.

{\bf Case 2: } $a  \not\in  V^{\bot},\: b  \not\in  V^{\bot},\:
 <a> + V^{\bot} \neq <b> + V^{\bot}$.

We first consider Case 1.  W.l.g, let $b \in V^{\bot}$. Then, $a$ is
not in $V^{\bot}$. So, Lemma \ref {PrpIncMat} implies

$\tilde{Q}_{\alpha, \beta} = s^{t-1}\sum \limits_{\gamma \in F}
Q_{\gamma,\beta}, \:\alpha, \beta \in F, Q = P,M.$

Now, by Theorems \ref {N-repl} and \ref {pairInc}  it follows that

$$ \tilde{n}^{ab}_{\alpha, \beta} =  s^{t-1} r^b_{\beta}
\mbox{ and }  \tilde{L_a} \tilde{L'_b}  = s^{t-1}k r^b_{\beta}.$$
Now, the result follows from Definition \ref {orthbl}.

In Case 2, applying (e) of Lemma \ref {PrpIncMat} we get $\tilde{Q}
= s^{t-2} c_q J, Q = P,M$, where $C_q$ is the total sum of the
entries of $Q$. But in view of Theorems \ref {N-repl} and \ref
{pairInc} $c_p = k c_m$. Hence the condition of Definition \ref
{orthbl} is satisfied.

 Proof of Case (b):(i) The result is immediate from (a) of
 Lemma \ref {PrpIncMat}.

 (ii) It is enough to show that when $ a  \not\in  V^{\bot},\: b  \not\in  V^{\bot},\:
 <a> + V^{\bot} <b> + V^{\bot}$ and $M$ is positive diagonal, then
every entry of $\tilde{M}$ is non-zero. But (d) of Lemma \ref
{PrpIncMat} says that for every $\alpha, \beta \in F$,  $
\tilde{m}_{\alpha, \beta}$ is $s^{t-1}$ times the sum of entries of
$M$ on a transversal. Since $M$ is positive diagonal, the result
follows.

\section{Plans for three-level factors}
     In this section we concentrate on three-level factors and
     assume that interactions involving three or more factors are not
     present. In the preceding section we have obtained the conditions on the
  developing procedure, so that the final plan satisfies desirable
  properties. However, in order that the size of the experiment is
  not too large, one needs also to have a good start. We have been
  able to find such a good start - a plan (${\cal P}$ )for a $3^4$ experiment on two
  blocks of size four each - as shown below.

 ${\cal P} = \begin{array}{ccccccccccc}
\mbox {Blocks} & |  && B_1 &  & |& &B_2 & & \\
 \hline      A  &| 0 & 1 & 1& 2&|&0 & 0 & 2 & 2\\
             B  &| 0 & 1 & 2& 0&|&2 & 1 & 1& 2\\
             C  &| 0 & 1 & 0& 1&|&1 & 2 & 0 & 2\\
             D  &| 0 & 0 & 1& 1&|&2 & 1 &2 & 0\\ \hline \end{array}$

{\bf Properties of ${\cal P}$ :} The effects satisfy the following
defining relations.

 $\mbox{Block} \equiv ABC \equiv AC^2D^2 \equiv AB^2D \equiv BC^2D.$

This implies that the main effects and the two-factor interactions
form the following alias classes.
\begin{eqnarray}\label{aliasCls}
  {\cal A}_1&=&\{A,B^2C^2,BD^2,CD\},\\
 {\cal A}_2 & = & \{B,A^2C^2,AD,CD^2\},\\
 {\cal A}_3&=&\{C,AD^2,A^2B^2,BD\},\\
 {\cal A}_4&=& \{D,AC^2,A^2B,B^2C\}. \end{eqnarray}

The relation among the classes are shown in the following graph,
where  adjacency represents orthogonality.

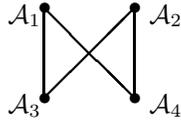
\begin{center}

$ \setlength{\unitlength}{4mm}
\begin{picture}(20,3)(0,2)

\thicklines

 \put(1,1){\line(0,1){3}} \put(1,1){\line(1,1){3}}

 \put(1,4){\line(1,-1){3}} \put(4,1){\line(0,1){3}}

\put(0.8,0.8){$\bullet$} \put(3.8,0.8){$\bullet$}
\put(0.8,3.8){$\bullet$} \put(3.8,3.8){$\bullet$}

\put(-0.2,0.5){$ {\cal A}_3$} \put(-0.2,3.5){${\cal A}_1$}
\put(4.5,0.5){${\cal A}_4$} \put(4.5,3.5){${\cal A}_2$}
 \end{picture} $

\end{center}

\vspace{.5em}

  We now modify ${\cal P}$ according to our requirement and expand
  along suitable subspaces to obtain the final plan $V({\cal P})$.
Before that we need a definition.

\begin{defi} Let us consider a plan $\rho$. Suppose  the set
${\cal E}$ of all effects of interest can be divided into
 several classes in such a way that every effect
is orthogonal to every other from a different class. Then $\rho$ is
called inter-class orthogonal and the classes will be referred to as
``orthogonal classes".\end{defi}

{\bf Remark 4.1:} The concept of  inter-class orthogonality has been
introduced and studied in Bagchi (1916).


  {\bf Case 1:} A $3^3$ experiment.
\begin{theo}\label{3^3} There exists an inter-class orthogonal plan for a $3^3$
experiment on six blocks of size four each estimating all main
effects and all  two-factor interactions.
\end{theo}

{\bf Proof :}
  We delete factor $D$ from ${\cal P}$
  to get a plan, say ${\cal P}^3$. Now the  alias classes are
  \begin{eqnarray*}
 {\cal A}_1&=&\{A,B^2C^2\},\\
  {\cal A}_2 & = & \{B,A^2C^2\},\\
  {\cal A}_3&=&\{C,A^2B^2\},\\
  {\cal A}_4&=& \{AC^2,A^2B,B^2C\}. \end{eqnarray*}

Now we expand ${\cal P}^3$ along $V = <\{(1,0,0)\}>$. Using Theorem
\ref {condOrthg} one can verify that in $V({\cal P}^3)$ the set of
all effects in the model satisfy inter-class orthogonality, the
orthogonal classes being as follows.
 \begin{eqnarray*}
C_1 = \{A, AC\}, &  C_2 = \{B, BC\}\\
C_3 = \{C, B^2C\} & \mbox{ and } C_4 = \{AB, AC^2, A^2B\}.
\end{eqnarray*}

Since no effect is aliased with any other, the result follows.
$\Box$

\vspace{.5em}

  {\bf Case 2:} A $3^4$ experiment.

\begin{theo}\label{3^4} There exists an inter-class orthogonal plan for a $3^4$
experiment on eighteen blocks of size four each estimating all main
effects and all two-factor interactions.
\end{theo}

{\bf Proof;} We start with ${\cal P}$ in which the alias classes are
as given in (\ref {aliasCls}) and the equations next to it. Now we
expand ${\cal P}$ along $V = <\{ (0,1,0,2), (1,0,1,0)\}>$. Theorem
\ref {condOrthg} implies that the resultant plan $V({\cal P})$ is
inter-class orthogonal for the set of all main effects and
two-factor interactions. The classes are :
  \begin{eqnarray*}
C_1&=&\{A,AC\}, \; C_2= \{B,BD^2\}\\
C_3&=&\{C \}, \; C_4 = \{D\}\\
C_5&=&\{BC,CD^2\}, \; C_6=\{AD,CD\}\\
C_7&=&\{AB,AD^2,\}, \; C_8=\{AB^2,BC^2\}\\
 C_{9} &=&\{AC^2,BD\}\Box
\end{eqnarray*}
Since no effect is aliased with any other, the result follows.
$\Box$ \vspace{.5em}

 {\bf Case 3:} A $3^5$ experiment.
\begin{theo}\label{3^5} There exists an inter-class orthogonal plan for a $3^5$
experiment on eighteen blocks of size four each estimating all main
effects and all two-factor interactions, except $DE^2$, which is
confounded with the block factor.
\end{theo}

{\bf Proof :} We obtain the following plan ${\cal P}^5$ from ${\cal
P}$ by adding a factor $E$.

  Plan  ${\cal P}^5 : \begin{array}{ccccccccccc}
\mbox {Blocks} & |  && B_1 &  & |& &B_2 & & \\
 \hline      A  &| 0 & 1 & 1& 2&|&0 & 0 & 2 & 2\\
             B  &| 0 & 1 & 2& 0&|&2 & 1 & 1& 2\\
             C  &| 0 & 1 & 0& 1&|&1 & 2 & 0 & 2\\
             D  &| 0 & 0 & 1& 1&|&2 & 1 &2 & 0\\
             E  &| 0 & 0 & 1& 1&|&2 & 1 &2 & 0 \\\hline\end{array}$

The effects satisfy the following defining relations :
 \begin{eqnarray*}
\mbox{Block} \equiv DE^2 \equiv ABC & \equiv AC^2D^2 \equiv AC^2E^2,
 \\ \equiv AB^2D \equiv AB^2E & \equiv BC^2D \equiv BC^2E^2.
\end{eqnarray*}

 Thus, interaction $DE^2$ is confounded with the block factor and the
 alias classes are as given below.
\begin{eqnarray*}
 {\cal A}_1&=&\{A,B^2C^2,BD^2,CD,CE,BE^2\},\\
 {\cal A}_2 & = & \{B,A^2C^2,AD,AE,CD^2,CE^2\},\\
 {\cal A}_3&=&\{C,AD^2,A^2B^2,AE^2,BD,BE\},\\
 {\cal A}_4&=& \{D,E,D^2E^2,AC^2,A^2B,B^2C\}. \end{eqnarray*}

We expand ${\cal P}$ along $V = <\{ (0,1,0,2,0), (1,0,1,0,2)\}>$ to
obtain the final plan $V({\cal P}^5)$. Using Theorem \ref
{condOrthg} we find that if we forget $DE^2$, then  $V({\cal P}^5)$
can be viewed as inter-class orthogonal plan : the classes being as
listed below.
  \begin{eqnarray*}
C_1&=&\{A,AC,CE^2\}, \; C_2= \{B,BD^2\}\\
C_3&=&\{C,E,AE^2\}, \; C_4 = \{D\}\\
C_5&=&\{BC,BE^2,CD^2\}, \; C_6=\{CD,AD\}\\
C_7&=&\{AD^2,AB,DE\}, \; C_8=\{BE,BC^2,AB^2\}\\
C_9&=&\{CE,AE\},\; C_{10}=\{BD,AC^2\}\Box
\end{eqnarray*}

{\bf Remark 4.2:} A resolution $V$ plan for a $3^4$ as well as a
$3^5$ experiment requires $81$ runs. Thus, apart from providing more
flexibility due to small size of the blocks, we have saved 9 runs in
both the situations. While all effects in the model are estimable in
the former plan, only one two-factor interaction is lost in the
later.

 {\bf Case 4:} A $3^6$ experiment.
\begin{theo} \label{3^6} There exists an inter-class orthogonal plan for a $3^5$
experiment on eighteen blocks of size four each estimating all main
effects and all but nine two-factor interactions. \end{theo}

{\bf Proof:}  We obtain the following plan ${\cal P}^6$ from ${\cal
P}$ by adding  factors $E$ and $F$..

  Plan  ${\cal P}^6 : \begin{array}{ccccccccccc}
\mbox {Blocks} & |  && B_1 &  & |& &B_2 & & \\
 \hline      A  &| 0 & 1 & 1& 2&|&0 & 0 & 2 & 2\\
             B  &| 0 & 1 & 2& 0&|&2 & 1 & 1& 2\\
             C  &| 0 & 1 & 0& 1&|&1 & 2 & 0 & 2\\
             D  &| 0 & 0 & 1& 1&|&2 & 1 &2 & 0\\
             E  &| 0 & 0 & 1& 1&|&2 & 1 &2 & 0\\
             F  &| 0 & 1 & 1& 2&|&0 & 0 & 2 & 2\\\hline\end{array}$

The effects satisfy the following defining relations :
 \begin{eqnarray*}
\mbox{Block} \equiv DE^2 \equiv AF^2 \equiv ABC  \equiv AB^2D \equiv
AB^2E \equiv AC^2D^2 \equiv AC^2E^2\\
 \equiv BC^2D \equiv BC^2E \equiv BCF \equiv B^2DF \equiv B^2EF \equiv C^2E^2F
\equiv C^2D^2F\end{eqnarray*}

  Thus, interactions
$DE^2$ and $AF^2$ are confounded with the block factor and the alias
classes are as given below.
\begin{eqnarray*}
 {\cal A}_1&=&\{A, A^2F^2,B^2C^2,BD^2,BE^2, CD,CE,F\},\\
 {\cal A}_2 & = & \{B,A^2C^2,AD,AE,CD^2,CE^2,C^2F^2,DF,EF\},\\
 {\cal A}_3&=&\{C,A^2B^2,AD^2,AE^2,BD,BE,B^2F^2,D^2F,E^F\},\\
 {\cal A}_4&=& \{D,E,A^2B,AC^2,B^2C,BF^2,C^2F,D^2E^2,\}. \end{eqnarray*}

We expand ${\cal P}^6$ along $V = <\{ (1,1,0,1,0,0)
(0,0,1,0,1,1)\}>$ to obtain the final plan $V({\cal P}^6)$. From
Theorem \ref {condOrthg} it follows that all the effects in the
assumed model, except $DE^2$ and $AF^2$ are divided into the
orthogonal classes listed below. The classes do contain pairs of
aliased two-factor interactions, which are presented within ().
\begin{eqnarray*}
C_1 = \{A,B,AD,\}, \; C_2 = \{C,E\}, \; C_3 = \{D,AB^2,BD\}, \\
C_4 = \{F,CE,CF,EF\},\; C_5 = \{BD^2,CF^2\}, C_6 = \{CD^2,BE^2\},\\
\; C_7 = \{BE, BF^2,D^2E^2\}, C_8 = \{(AD^2,EF^2) , (AB^2, CF^2)\}, \;  \\
C_9 = \{AE^2, BC^2, (AC^2,BF^2)\}, \;C_{10} = \{CD, AC, (AE, DF), (AF,BC)\}, \; \\
\end{eqnarray*}
We see that there are five pairs of mutually aliased two-factor
interactions. These together with the interactions confounded with
the blocks make the size of the set on non-estimable two-factor
interactions as seven. Again, the classes $C_4$ and $C_{10}$ are of
size four (counting each alias pair as one effect). But in the
present set up at most three non-orthogonal effects can be
estimated, so that two more two-factor interactions are lost. Hence
the result. $\Box$

{\bf Remark 4.3:} In a plan of 18 blocks of size 4 each, the
available treatment degrees of freedom is 54 and so at most 27
effects can be estimated. In $V({\cal P}^6)$ these  treatment
degrees of freedom are utilized to the full extent as 6 main effects
and 30-9 = 21 two-factor interactions are estimated.

  We now present another plan, which
  may be viewed as a supplement of the plan $V({\cal P}^6)$ of
  Theorem \ref {3^6}, in the sense that these two plans together
  estimates all main effects and all two-factor interactions.

  \begin{theo} \label{N3^6}There exists a plan $V({\cal P}_2^6)$ on 6 blocks of size 4 each,
  which  estimates  all but one of the two-factor interactions lost in
   $V({\cal P}^6)$, so that these two plans together estimates all
   effects in the model except one  two-factor interaction.
\end{theo}

{\bf Proof :} The plan is as given below.

$\begin{array}{ccccccccccc}
\mbox {Blocks} & |  && B_1 &  & |& &B_2 & & \\
 \hline      A  &| 0 & 1 & 1& 2&|&0 & 0 & 2 & 2\\
             B  &| 0 & 1 & 2& 0&|&2 & 1 & 1& 2\\
             C & |0 & 1 & 1& 2&|&0 & 0 & 2 & 2\\
             D  &| 0 & 1 & 2& 0&|&2 & 1 & 1& 2 \\
             E  &| 0 & 0 & 1& 1&|&2 & 1 &2 & 0\\
             F  &| 0 & 1 & 0& 1&|&1 & 2 & 0 & 2 \\\hline\end{array}$

We can see that the effects $AC^2$ and $BD^2$ are confounded with
the block factor. The remaining effects form the following alias
classes.
\begin{eqnarray*}
 {\cal A}_1&=&\{A,C, AC,BE^2,B^2F^2,  DE^2, D^2F^2,   EF\},\\
 {\cal A}_2 & = & \{B, D, AE,  A^2F^2,BD,CE,C^2F^2,E^2F\},\\
 {\cal A}_3&=&\{F, A^2B^2,A^2D^2,AE^2,B^2C^2, BE,C^2D^2, CE^2, DE\},\\
 {\cal A}_4&=& \{E,A^2B, A^2D, AF^2, BC^2,B^2F, C^2D, CF^2, D^2F,\}. \end{eqnarray*}

Now, we expand the above plan along $V = <\{(1,0,0,1,2,0)\}>$ and
get the required plan $V({\cal P}_2^6)$. We now present the
orthogonal classes obtained in view of Theorem \ref {condOrthg}. The
effects within () are aliased.
\begin{eqnarray} \label{orthCls} C_1 = \{(C,DE^2,BF), (B,AE,CF)\},\\
~~~~C_2 = \{(A,AC,BE^2), (D,BD,EF^2),(EF,DF),(AF,CE)\},\\
C_3 =\{(F,DE,BC),(AD^2,BC^2,BF^2,CF^2)\},\\
C_4 = \{(E,AB^2,DF^2), (AD,CE^2),(CD^2,AF^2),
\\(AB,AE^2,BE,CD)\}\end{eqnarray}

 We shall now see how the lost effects can be estimated. When two
 effects are aliased in  $V({\cal P}^6)$ we need to see that one of them can
 be  estimated from $V({\cal P}_2^6)$. When there are four effects
 among which at most three can be estimated from $V({\cal P}^6)$, we
 need to find one member which can be  estimated from $V({\cal P}_2^6)$. We
 present the details in the following table.

\begin{center}
 $\left[\begin{array}{cccccccc} \hline
\mbox{Effects to}     &\mbox{ Class of the } &\mbox{Earlier} \\
\mbox{be Estimated}   &\mbox{New plan} &\mbox{Confounded
with}\\\hline
   DE^2     & C_1 & \mbox{Block} \\
   AF^2     & C_4 & \mbox{Block} \\ \hline
   EF^2     & C_2 & AD^2 \\
   DF       & C_2  & AE    \\
   BF^2     & C_3  & AC^2       \\
   BC       & C_3  & AF    \\
   AB^2     & C_4  & CF^2     \\\hline
  CF       & C_1  & \mbox{one of } F, CE, EF   \\
   CD       & C_4  & \mbox{one of } AC,BC,DF    \\
                \hline
\end{array} \right ].$ \end{center}

We note that $C_4$ contains three and  each of the other classes
contains two of the effects to be estimated. Further, one may check
from equation (\ref {orthCls} ) and the following three equations
that the effects in the table, which belongs to the same class of
$V({\cal P}_2^6)$ are not aliased. But since from every class at
most two effects can be estimated, one effect in $C_4$ is lost.
Thus, all but one of the lost nine effects can be estimated. Hence
the result. $\Box$

{\bf Remark 4.4} It is known that a resolution $V$ plan for a $3^6$
experiment requires $3^5 = 243$ runs [see Theorem 13.1 of Hinkelman
and Kempthorn (205), for instance]. The plans $V({\cal P}^6)$ and
$V({\cal P}_2^6)$ together provides estimability almost to the same
extent, but requires only 96 runs.

  {\bf Acknowledgement :}

\section{References}

\begin{enumerate}
\item Addelman,S. (1962). Orthogonal main effect plans for asymmetrical
factorial experiments. Technometrics, vol. 4, p: 21-46.

\item  Bagchi, S. (2010). Main effect plans orthogonal through the
block factor, Technometrics, vol. 52, p : 243-249.

\item Bagchi, S. (2016) Inter-class orthogonal main effect plans for
asymmetrical experiments. To appear in Sankhya.

\item Bose, M. and  Bagchi, S. (2007). Optimal main effect plans on blocks
of small size, Jour. Stat. Prob. Let., vol.  77, p : 142-147.

\item Bose, R.C. (1947).  Mathematical theory of the
symmetrical factorial design. Sankhya, vol. 8, p : 107-166.



\item Hinkelman, K. and Kempthorn, O. (2005). Design and analysis of
experiments, vol.2, Wiley Series in Probability and Statistics.

\item Morgan, J.P. and Uddin, N. (1996). Optimal blocked main effect
plans with nested rows and columns and related designs, Ann. Stat.,
vol. 24, p : 1185-1208.


\item Mukerjee, R., Dey, A. and Chatterjee, K. (2001). Optimal
main effect plans with non-orthogonal blocking.  Biometrika, vol.
89, p : 225-229.

\end{enumerate}

\end{document}